\documentclass{article}
 \usepackage{amssymb}
\usepackage{amsmath}
\title{Boundedness of the gradient of a solution  to the Neumann-Laplace problem in a convex domain}
\author{Vladimir Maz'ya\footnote{The author was partially supported by  the UK Engineering and Physical
Sciences Research Council grant EP/F005563/1.}
\\*[10pt]
\emph{\small Department of Mathematical Sciences,   University of Liverpool, Liverpool L69 7ZL}
\\
\emph{\small and}\\
\emph{\small Department of Mathematics, Link\"oping University, Link\"oping, SE-581 83}
\\
\emph{\small e-mail: vlmaz@mai.liu.se}
}
\date{}

\begin{document}
\maketitle

\begingroup

\narrower\noindent
{\it Abstract.}  It is shown that solutions of the Neumann problem for the Poisson equation in an arbitrary convex $n$-dimensional domain  are uniformly Lipschitz.  
Applications of this result to some aspects of regularity of solutions to the Neumann problem on convex polyhedra are given.
\endgroup


\medskip\noindent
{ {\small\it Key words}: }{\small Neumann problem, Laplace operator, Sobolev spaces, convex domain, boundedness of the gradient, eigenvalues of the    Neumann-Beltrami  operator

\medskip

\begingroup
\narrower\noindent
{\it R\'esum\'e.}
On d\'emontre que les solutions du probl\`eme de Neumann pour l'\'equation
de Poisson
dans un domaine convexe arbitraire de dimension $n$ sont uniform\'ement
Lipschitz.
Les applications de ce r\'esultat \`a quelques aspects de r\'egularit\'e de solutions du probl\`eme de Neumann sur les poly\`edres convexes sont donn\'ees.
\endgroup


\medskip\noindent
{ {\small\it Mots cl\'e}}: probl\`eme de Neumann, op\'erateur de Laplace, espaces de Sobolev,
domaine convexe; bornitude du gradient, valeurs propres de l'op\'erateur
de Neumann-Beltrami.
\section{Introduction}

Let $\Omega$ be a bounded convex domain in $\Bbb{R}^n$ and let $W^{l,p}(\Omega)$ stand for the Sobolev space of functions in $L^p(\Omega)$ with distributional derivatives of order $l$ in $L^p(\Omega)$.  By $L_\bot^p(\Omega)$ and $W_\bot^{l,p}(\Omega)$ we denote the subspaces of functions $v$ in $L^p(\Omega)$  and $W^{l,p}(\Omega)$ subject to $\int_\Omega v dx =0$.

\smallskip

Let $f\in L^2_\bot(\Omega)$   and let $u$ be  the unique function in $W^{1,2}(\Omega)$, also  orthogonal to $1$ in 
$L^2(\Omega)$,  and satisfying the Neumann problem
\begin{eqnarray}\label{s1}
-\Delta u &=& f \quad {\rm in} \,\,\, \Omega,\\
\frac{\partial u}{\partial \nu} & = & 0 \quad {\rm on}\,\,\, \partial\Omega,\label{s2}
\end{eqnarray}
where $\nu$ is the unit outward normal vector to $\partial \Omega$ and the problem (\ref{s1}), (\ref{s2})  is understood  in the variational sense. It is well known that the inverse mapping
\begin{equation}\label{s3}
L^2_\bot(\Omega) \ni f\to u \in W_\bot^{2,2}(\Omega)
\end{equation} 
is continuous. (Since any attempt at reviewing the rich history of this fact and other ones, closely related to it, within  frames of a short article is hopeless, I restrict myself to a number of relevant references \cite{Be}, \cite{BS}, \cite{Gr}, \cite{Ka}, \cite{L1}--\cite{LU}, \cite{M67}, \cite{M73}, \cite{Sch}, \cite{So}.)  As shown in \cite{AJ} (see also \cite{JMM} for a different proof, and \cite{Ad}, \cite{Fr1}--\cite{FJ}  for the Dirichlet problem), the operator
\begin{equation}\label{s3a}
L^p_\bot(\Omega)  \ni f\to u \in W_\bot^{2,p}(\Omega)
\end{equation} 
is also continuous if $1<p<2$. One cannot guarantee the continuity of (\ref{s3a}) for any $p\in (2,\infty)$ without additional information about the domain. The situation is the same as in the case of the Dirichlet problem (see \cite{BS}, \cite{Fr1}-\cite{FJ}), which, moreover, possesses the following useful property: if $\Omega$ is convex, the gradient of the solution is uniformly bounded provided the right-hand side of the equation is good enough. This property can be easily checked by using a simple barrier. Other approaches to similar results were exploited in \cite{M69} and \cite{KM}  for different equations and systems but only for the Dirichlet boundary conditions.

\smallskip

In this respect, other boundary value problems are in a nonsatisfactory state. For instance,  it was unknown up to now whether solutions of the problem (\ref{s1}), (\ref{s2}) with a smooth $f$ are uniformly Lipschitz under the only condition of convexity of $\Omega$. 

\smallskip

The main result of the present paper is {\it the boundedness of $|\nabla u|$ for the solution $u$ of the Neumann problem $(\ref{s1}), (\ref{s2})$ in any convex domain} $\Omega \subset \Bbb{R}^n$, $n\geq 2$. 

\smallskip

 A  direct consequence of this fact is the sharp lower estimate $\Lambda\geq n-1$ for the first nonzero eigenvalue $\Lambda$ of the Neumann problem for the  Beltrami  operator on a convex subdomain of a unit hemisphere. It was  obtained by a different argument  for manifolds of positive Ricci curvature by J. F. Escobar in \cite{Esc}, where the case of equality was settled as well.  This estimate  answered a question  raised by M. Dauge \cite{Da},  and it leads, in combination with known techniques of the theory of elliptic boundary value problems in domains with piecewise smooth boundaries (see \cite{Da}, \cite{MP1}--\cite{MR}), to estimates for solutions of the problem (\ref{s1}), (\ref{s2}) in various function spaces. Two examples are given at the end of this article.

\medskip

\section{Main result}

In what follows,  we need a constant $C_\Omega$ in the relative isoperimetric inequality
\begin{equation}\label{s6}
s(\Omega\cap \partial g) \geq C_\Omega\, |g|^{1-1/n},
\end{equation} 
where $g$ is an arbitrary open set in $\Omega$ such that $|g|\leq |\Omega|/2$ and $\Omega\cap \partial g$ is a smooth (not necessarily compact) submanifold of $\Omega$. By $s$ we denote the $(n-1)$-dimensional area and by $|g|$ the $n$-dimensional Lebesgue measure. The Poincar\' e-Gagliardo-Nirenberg inequality
\begin{equation}\label{s7}
\inf\limits_{t\in \Bbb{R}} \|v-t\|_{L^{n/(n-1)}(\Omega)} \leq  const. \, \|\nabla v \|_{L^1(\Omega)}, \qquad \forall v\in W^{1,1}(\Omega),
\end{equation} 
where $const. \leq C^{-1}_\Omega$ is a consequence of (\ref{s6}) (see Theorem 3.2.3 \cite{M85}). 


\smallskip

{\bf Theorem.} {\it Let $f\in L_\bot^q(\Omega)$ with a certain $q>n$. Then the solution $u\in W_\bot^{1,2}(\Omega)$ of the problem $(\ref{s1}), (\ref{s2}) $ satisfies the estimate
\begin{equation}\label{s9}
\|\nabla u\|_{L^\infty(\Omega)} \leq c(n,q) \,  C^{-1}_\Omega \, |\Omega|^{(q-n)/qn} \|f\|_{L^q(\Omega)}, 
\end{equation} 
where $c$ is a constant depending only on $n$ and $q$. }

{\bf Proof.} It suffices to prove (\ref{s9}) assuming that $f\in C_0^\infty(\Omega)$. Let us approximate $\Omega$ by a sequence $\{\Omega_m\}_{m\geq 1}$ of convex domains with smooth boundaries, $\Omega_m\supset \overline{\Omega}$.  This can be done, for instance, by approximating $\Omega$ by a family of equidistant surfaces and by smoothing them with small perturbation of normal vectors. Then
(\ref{s7}) implies
\begin{equation}\label{s10}
\inf\limits_{t \in \Bbb{R}} \|w-t \|_{L^{\frac{n}{n-1}}(\Omega_m)} \leq (1+\varepsilon)\, C^{-1}_\Omega \, \|\nabla w \|_{L^1(\Omega_m)},
\end{equation} 
for all $w\in W^{1,1}(\Omega_m)$, where $\varepsilon$ is an arbitrary positive number and $m= m(\varepsilon)$.

\smallskip

By $u_m$ we denote a solution of the problem (\ref{s1}), (\ref{s2}) in $\Omega_m$ with $f$ extended by zero outside $\Omega$. One can easily see that $\nabla u_m \to \nabla u$ in $L^2(\Omega)$. Hence, it is enough  to obtain (\ref{s9}) assuming that $\partial\Omega$ is smooth. 

\smallskip

Let $t>\tau>0$ and let $\Psi$ be a piecewise linear continuous function on $\Bbb{R}$ specified by $\Psi (\xi) = 0$ for $\xi<\tau$ and $\Psi (\xi) = 1$ for $\xi> t$. Note that
\begin{equation}\label{s11}
(\Delta u)^2 - |\nabla_2 u|^2 = ( u_{x_j} \, \Delta u)_{x_j} - ( u_{x_j} u_{x_i x_j})_{x_i},
\end{equation} 
where
$$|\nabla_2 u| = \Bigl( \sum_{1\leq i,j\leq n} u_{x_ix_j}^2 \Bigr)^{1/2}.
$$
Hence
\begin{eqnarray}\label{s12}
&& \int_\Omega \Psi(|\nabla u|) ( f^2 - |\nabla_2 u|^2) \, dx = \int_{\partial\Omega} \Psi(|\nabla u|)\bigl(\nu_j \,  u_{x_j} \Delta u - \nu_i \, u_{x_j}\,  u_{x_i x_j}\bigr) \, ds_x \nonumber\\
&&+\int_\Omega \Psi'(|\nabla u|) \bigl ( \, (|\nabla u| )_{x_j} u_{x_j} f + (|\nabla u| )_{x_i}\, u_{x_j}\, u_{x_ix_j} \bigr) dx,
\end{eqnarray}
where $(\nu_1, \ldots , \nu_n)$ are components of the outward unit normal. By the Bernshtein-type identity (see, for instance, \cite{Gr} or \cite{L2}), the first integral on the right-hand side of (\ref{s12}) equals
$$-2\int_{\partial\Omega} Q(\nabla_{\it tan} u, \, \nabla _{\it tan} u) \, ds_x,$$
where $Q$  is the second fundamental quadratic form on $\partial\Omega$ and $ \nabla _{\it tan}$ is the tangential gradient. The form $Q$ is nonpositive by convexity of $\Omega$, which leads, together with (\ref{s12}), to the inequality
\begin{equation}\label{s13}
\int_\Omega \Psi'(|\nabla u|) \bigl ( \, (|\nabla u| )_{x_j} u_{x_j} f + (|\nabla u| )_{x_i}\, u_{x_j}\, u_{x_ix_j} \bigr) dx \leq \int_\Omega \Psi(|\nabla u|)\, f^2\, dx.
\end{equation} 

\smallskip

By the co-area formula \cite{Fe}, the left-hand side of (\ref{s13}) is identical to
$$(t-\tau)^{-1}\int_\tau^t \int_{|\nabla u|=\sigma} \bigl ( \, (|\nabla u| )_{x_j} u_{x_j} f + (|\nabla u| )_{x_i}\, u_{x_j}\, u_{x_ix_j} \bigr)  \frac{ds_x}{|\nabla|\nabla u|\, |} d\sigma,$$
which is equal to
\begin{equation}\label{s14}
(t-\tau)^{-1}\int_\tau^t \int_{|\nabla u|=\sigma} \Bigl ( \frac{\partial u}{\partial \nu} f - |\nabla u|\frac{\partial}{\partial \nu} |\nabla u| \Bigr) ds_x\, d\sigma
\end{equation} 
because $\nabla |\nabla u| = -\nu\,  |\nabla|\nabla u|\, |$ on the level surface $|\nabla u| =\sigma$, where $\nu$ is the unit normal, outward with respect to the set $\{x: |\nabla u|>\sigma\}$. The expression (\ref{s14}) can be written as
$$(\tau -t)^{-1} \int_\tau^t \int_{|\nabla u|=\sigma} \Bigl ( \frac{\partial u}{\partial \nu} f - |\nabla u| \, |\nabla |\nabla u|\, | \Bigr) ds_x\, d\sigma.$$
Passing here to the limit as $\tau\uparrow t$ and using (\ref{s13}), we arrive at the estimate
$$\int_{|\nabla u|=t}  \Bigl( |\nabla u| \, |\nabla |\nabla u|\, |  - f \frac{\partial u}{\partial\nu} \Bigr) ds_x \leq \int_{|\nabla u|>t} f^2 dx,$$
which implies
\begin{equation}\label{s15}
t\int_{|\, \nabla u|=t}   |\nabla |\nabla u|\, | ds_x \leq t \int_{|\nabla u|=t} |f|\, ds_x + \int_{|\nabla u|>t} f^2 dx.
\end{equation} 
 We define the median of $|\nabla u|$ as
 $${\it med}\,  |\nabla u| = \sup \{t\in \Bbb{R} : |\{|\nabla u| >t \} | \geq |\Omega|/2\}, $$
 and we note that
 $$ \bigl |\{|\nabla u  | >{\it med}\,  |\nabla u|  \} \bigr | \leq |\Omega|/2$$
and
$$ \bigl |\{|\nabla u  |  \geq {\it med}\,  |\nabla u|  \} \bigr | \geq |\Omega|/2.$$

\noindent
 Clearly,
 \begin{equation}\label{s16}
{\it med} \, |\nabla u|  \leq \Bigl(\frac{2}{|\Omega|}\Bigr)^{1/2} \, \|\nabla u\|_{L^2(\Omega)}
\end{equation} 
and, by H\"older's inequality and (\ref{s10}), 
\begin{eqnarray*}
\|\nabla u\|_{L^2(\Omega)} & \leq &  \inf\limits_{\gamma\in \Bbb{R}} \|u-\gamma\|_{L^{\frac{n}{n-1}}(\Omega)}^{1/2} \|f\|_{L^n(\Omega)}^{1/2}\\
&\leq &  (1+\varepsilon)\, C_\Omega^{-1/2}  \|\nabla u\|_{L^1(\Omega)}^{1/2}  \|f\|_{L^n(\Omega)}^{1/2}.
\end{eqnarray*}
Hence,
$$\|\nabla u\|_{L^2(\Omega)}  \leq (1+\varepsilon)^2  C_\Omega^{-1} |\Omega|^{1/2}\,  \|f\|_{L^n(\Omega)},$$

\smallskip

\noindent
which, in combination with (\ref{s16}) and H\"older's inequality, implies
\begin{equation}\label{s17}
{\it med}\,  |\nabla u|  \leq 2^{1/2}(1+\varepsilon)^2 C_\Omega^{-1} |\Omega|^{(q-n)/qn} \, \|f\|_{L^q(\Omega)}.
 \end{equation}
 We introduce the function
 $$\psi : \bigl [{\it med} \, |\nabla u|, \, \max |\nabla u| \,\bigr ] \to [ 0, \infty)$$
 by the equality
 \begin{equation}\label{s18}
 \psi(t) = \int_{{\it med} \, |\nabla u|}^t \Bigl( \int_{|\nabla u|=\sigma} |\nabla |\nabla u|\, | ds_x\Bigr)^{-1} d\sigma.
  \end{equation}

\noindent
 Let  $E_\psi = \bigl\{ x: |\nabla u(x)| = t(\psi)\bigr\}$ and $M_\psi = \bigl\{ x: |\nabla u(x)| > t(\psi)\bigr\}$. 
  Putting $t= t(\psi)$ in (\ref{s15}) and 
 integrating  in $\psi$ over $\Bbb{R}_+$, we arrive at 
 $$\max |\nabla u|^2 \leq ({\it med}\,  |\nabla u|)^2 + 2\max |\nabla u| \int_0^\infty \int_{E_\psi} |f|\, ds_x\, d\psi + 2 \int_0^\infty \int_{M_\psi} f^2\, dx\, d\psi.$$
Recalling (\ref{s17}), we see that in order to obtain (\ref{s9}), it suffices to prove the inequalities 
  \begin{equation}\label{s19}
 \int_0^\infty \int_{E_\psi} |f|\, ds_x\, d\psi \leq c\, C_\Omega^{-1} |\Omega|^{(q-n)/qn}  \|f\|_{L^q(\Omega)}
 \end{equation}
 and 

 \begin{equation}\label{s20}
\Bigl(\int_0^\infty \int_{M_\psi} f^2\, dx\, d\psi \Bigr)^{1/2} \leq c\, C_\Omega^{-1} |\Omega|^{(q-n)/qn}  \|f\|_{L^q(\Omega)},
 \end{equation}
where $c= c(n,q)$.

\smallskip
 
 The following argument leading to (\ref{s19}) is an obvious modification of the proof of Lemma 4 \cite{M69a}. We start with the estimate for the area of the set $E_\psi$
 \begin{equation}\label{s21}
 s(E_\psi)^2 \leq -\frac{d}{d \psi} |M_\psi |
  \end{equation}
  obtained in Lemma 2 \cite{M69a}. By the triple H\"older inequality and (\ref{s21})
 $$
 \int_{E_\psi} |f|\, ds_x \leq \Bigl( \int_{E_\psi} |f|^q \frac{ds_x}{|\nabla |\nabla u|\, |}\Bigr)^{1/q} \Bigl(\int_{E_\psi} |\nabla |\nabla u|\, |\, ds_x\Bigr)^{1/q} s(E_\psi)^{1-2/q}$$
$$\leq  s(E_\psi)^{-1} \Bigl( -\frac{d}{d\psi} |M_\psi|/ t'(\psi)\Bigr)^{1-1/q} \Bigl( \int_{E_\psi} |f|^q \frac{ds_x}{|\nabla |\nabla u|\, |}\Bigr)^{1/q}\int_{E_\psi} |\nabla |\nabla u|\, |\, ds_x.$$
This, in combination with the inequality
$$ 
s(E_\psi) \geq (1+\varepsilon)^{-1} C_\Omega \, |M_\psi|^{1-1/n}$$
(see (\ref{s6})), implies the estimate for the integral on the left-hand side of (\ref{s19}) 
 $$ \int_0^\infty \int_{E_\psi} |f|\, ds_x\, d\psi \leq  (1+\varepsilon)C_\Omega^{-1}\Bigl(\int_0^{|\Omega|/2} \mu^{\frac{(1-n)q}{n(q-1)}} d\mu\Bigr)^{1-1/q} \Bigl(\int_0^\infty \int_{E_\psi} |f|^q \frac{ds_x}{|\nabla |\nabla u|\, |}\, t'(\psi)\, d\psi\Bigr)^{1/q}$$
  $$\leq (1+\varepsilon) C_\Omega^{-1}\Bigl(\frac{n(q-1)}{q-n} \Bigl(\frac{|\Omega|}{2}\Bigr) ^{\frac{q-n}{n(q-1)}}\Bigr)^{1-1/q} \|f\|_{L^q(\Omega)}$$
  and the proof of (\ref{s19}) is complete. 
  
  \smallskip
  
  We turn to inequality (\ref{s20}). By (\ref{s21}), its left-hand side does not exceed
  $$\Bigl(-  \int_0^\infty \int_{M_\psi} |f|^2 dx \frac{d\, |M_\psi|}{s(E_\psi)^2} \Bigr)^{1/2}$$
  which  is dominated by
  $$(1+\varepsilon)\, C_\Omega^{-1}\Bigl(-  \int_0^\infty \int_{M_\psi} |f|^2 dx \frac{d\, |M_\psi|}{|M_\psi|^{2\frac{n-1}{n}}} \Bigr)^{1/2}$$
  $$ \leq (1+\varepsilon)\, C_\Omega^{-1}\Bigl(\int_0^{|\Omega|/2} \int_0^\sigma f_*(\tau)^2 d\tau \frac{d\sigma}{\sigma^{2\frac{n-1}{n}}} \Bigr)^{1/2},$$
 where $f_*$ is the nonincreasing rearrangement of $|f|$. Now, (\ref{s20}) follows by integration by parts and H\"older's inequality. The proof of  Theorem is complete. 
 
\medskip 

\section{Regularity of solutions to the Neumann problem in a convex polyhedron}
The following assertion essentially stemming from the above theorem is a particular case of Escobar's result in \cite{Esc} mentioned in Introduction. 

{\bf Corollary.} {\it Let $\omega$ be a convex subdomain of the upper unit hemisphere $S_+^{n-1}$. The first positive eigenvalue $\Lambda$ of the Beltrami operator on $\omega$ with zero Neumann data on $\partial\omega$ is not less than} $n-1$.

{\bf Proof.} Let $\lambda(\lambda + n-2) = \Lambda$ and $\lambda >0$. In the convex domain
$$\Omega =\bigl \{x\in \Bbb{R}^n: 0< |x| <1, \, \frac{x}{|x]} \in \omega \bigr\},$$
we define the function
\begin{equation}\label{s90}
u(x) = |x|^\lambda \Phi\Bigl (\frac{x}{|x|}\Bigr) \eta (|x|),
\end{equation}
where $\Phi$ is an eigenfunction corresponding to $\Lambda$ and $\eta$ is a smooth cut-off function on $[0, \infty)$, equal to one on $[0, 1/2]$ and vanishing outside $[0, \, 1]$. 

\smallskip

Let $N$ be an integer satisfying $4N > n-1 \geq 4(N-1)$ and let $j= 0, 1, \ldots, N$,
$$q_j =
\begin{cases}
\displaystyle{\frac{2(n-1)}{n-1-4j}}, & {\rm if} \,\, 0\leq j<  (n-1)/4,\\
{\rm arbitrary} & {\rm if}\,\,   j= (n-1)/4,
\end{cases}
$$
and $q_N =\infty$.  Iterating the estimate 
$$\|\Phi \|_{L^{q_j +1}(\omega)} \leq c\,  \Lambda\,  \|\Phi \|_{L^{q_j}(\omega)}$$
obtained in Theorems 5 and 6 \cite{M69a}, we see that  $\Phi \in L^\infty(\omega).$

\smallskip

The function $u$, defined by (\ref{s90}), satisfies the Neumann problem  (\ref{s1}), (\ref{s2}) with
$$f(x) = - \Phi\Bigl (\frac{x}{|x|}\Bigr) \bigl [ \Delta, \eta(|x|) \bigr] \, |x|^\lambda.$$
Since $\Phi \in L^\infty (\omega)$, it follows that $f\in L^\infty(\Omega)$ and by Theorem, $|\nabla u| \in L^\infty(\Omega)$, which is possible only if $\lambda\geq 1$, i.e. $\Lambda  \geq n-1$.  The proof is complete. 

\medskip



 Two applications of the above estimate for $\Lambda$
will be formulated. 

\smallskip

Let $\Omega$ be a convex bounded $3$-dimensional polyhedron. By the techniques, well-known nowadays (see \cite{Da},  \cite{MP1}--\cite{MR}), one can show   the unique solvability of the variational Neumann problem in $W_\bot^{1,p}(\Omega)$ for every $p\in (1, \infty)$.  By definition of this problem, its solution is subject to the integral identity
$$\int_\Omega \nabla u\cdot \nabla \eta\, dx = f(\eta),$$
where $f$ is a given distribution in the  space $(W^{1, p'}(\Omega) )^*$, $ f(1) =0$ and $\eta$ is an arbitrary function in $W^{1, p'}(\Omega)$, $p+p' = p p'$. 

\smallskip

Let us turn to the second application of Corollary. We continue to deal with the polyhedron $\Omega$ in $\Bbb{R}^3$. Let $\{{\cal O}\}$ be  the collection of all vertices and let $\{U_{{\cal O}}\}$ be an open  finite covering of $\overline{\Omega}$ such that ${\cal O}$ is the only vertex in $U_{{\cal O}}$. Let also $\{E\} $ be the collection of all edges and let $\alpha_E$ denote the opening of the dihedral angle with edge $E$, $0<\alpha_E <\pi$. The notation $r_E (x)$ stands for the distance between $x\in U_{{\cal O}}$ and the edge $E$ such that ${\cal O} \in \overline{E}$. 

\smallskip

With every vertex ${\cal O}$ and edge $E$ we  associate real numbers $\beta_{\cal O}$ and $\delta_E$, and we introduce the weighted $L^p$-norm
$$\|v \|_{L^p(\Omega; \{\beta_{\cal O}\}, \{\delta_E\})} : = \Bigl( \sum _{\{{\cal O}\}} \int_{U_{{\cal O}}} |x - {\cal O}|^{p\beta_{\cal O}} \prod_{\{E: {\cal O} \in \overline{E}\}} \Bigl(\frac{r_E (x)}{|x - {\cal O}|} \Bigr)^{p\delta_E} |v (x)|^p dx\Bigr)^{1/p},$$
where $1<p<\infty$. Under the conditions
\begin{eqnarray*}
3/p' & > &  \beta_{\cal O} > -2 + 3/p',\\
2/p' & > &\delta_E > -\min\{2, \, \pi/\alpha_E\} + 2/p',
\end{eqnarray*}
the inclusion of the function $f$ in $L^p(\Omega; \{\beta_{\cal O}\}, \{\delta_E\})$ implies the unique solvability of the problem (\ref{s1}), (\ref{s2}) in the class of functions with all derivatives of the second order belonging to $L^p(\Omega; \{\beta_{\cal O}\}, \{\delta_E\})$. This fact follows from Corollary and a result essentially  established in Sect. 7.5 \cite{MR}.

\smallskip

An important particular case when all $\beta_{\cal O}$ and $\delta_E$ vanish, i.e. when we deal with a standard Sobolev space $W^{2,p}(\Omega)$, is also included here. To be more precise, if
\begin{equation}\label{s30}
1<p<\min\Bigl\{ 3, \frac{2\alpha_E}{(2\alpha_E -\pi)_+} \Bigr\}
\end{equation}
for all edges $E$, then the inverse operator of the problem  (\ref{s1}), (\ref{s2}):
$$
L_\bot^p(\Omega)\ni f \to u\in W_\bot^{2,p} (\Omega)$$
is continuous whatever the convex polyhedron $\Omega\subset \Bbb{R}^3$ may be.  The bounds for $p$ in (\ref{s30}) are sharp for the class of all convex polyhedra.

\end{document}